\input amssymb.sty
\documentclass{amsproc}
\usepackage{graphicx}

\theoremstyle{definition}

\theoremstyle{remark}

\numberwithin{equation}{section}

\begin{document}

\title{Nombre de factorisations d'un grand cycle}

\author{Philippe Biane}
\address{CNRS, D\'epartement de Math\'ematiques et Applications,
 \'Ecole Normale Sup\'erieure, 45, rue d'Ulm 75005 paris, FRANCE
}
\email{Philippe.Biane@ens.fr}
\subjclass{Primary ; Secondary }
\date{}

\begin{abstract}
On donne une d\'emonstration simple d'une formule de Goupil et Schaeffer
qui compte  le nombre de factorisations d'un cycle de longueur maximale dans
$S_n$ en produit de deux permutations de classes de conjugaisons donn\'ees.
\end{abstract}

\maketitle
Dans la suite on utilise les notations du livre 
de MacDonald \cite{M}.

Soit $c^n_{\lambda\mu}$ le nombre de factorisations dans $S_n$ 
d'un cycle de longueur $n$
en un produit de deux permutations de classes de conjugaison $\lambda$ et $\mu$.
 La th\'eorie des caract\`eres donne la formule 
 \begin{equation}\label{00}c^{\nu}_{\lambda\mu}=\frac{n!}{
z_{\lambda}z_{\mu}}\sum_{\rho\vdash n}\frac{
\chi_{\lambda}^{\rho}\chi_{\mu}^{\rho}\chi_{\nu}^{\rho}}{\chi^{\rho}_{1^n}}.
\end{equation}
pour le nombre de d\'ecompositions d'une permutation de classe $\nu$ en produit
de deux permutations de classes $\lambda$ et $\nu$. La somme porte sur les
partitions de $n$ et les $\chi^{\rho}$ sont les caract\`eres du groupe
sym\'etrique, tandis que $z_{\lambda}=\prod_i\alpha_i!\,i^{\alpha_i}$
si $\lambda=1^{\alpha_1}2^{\alpha_2}\ldots n^{\alpha_n}$.

Lorsque $\nu=(n)$, un cycle de longueur maximale, on a $\chi_\nu^{\rho}=0$ sauf
si $\rho$ est une \'equerre, c'est-\`a-dire un diagramme de la forme
$1^r(n-r)$, et on obtient
\begin{equation}\label{01}
 c^n_{\lambda\mu}=\frac{n}{
z_{\lambda}z_{\mu}}\sum_{r=0}^{n-1}(-1)^rr!(n-1-r)!
\chi_{\lambda}^{1^r(n-r)}\chi_{\mu}^{1^r(n-r)}.
\end{equation}
qui est la formule (4) de \cite{GS}. Cet article se poursuit par une analyse 
combinatoire de cette formule, pour la transformer en une expression ne
contenant
que des termes positifs.
Nous allons suivre une voie plus alg\'ebrique et
introduire une fonction g\'en\'eratrice pour ces quantit\'es en
utilisant les fonctions sym\'etriques $p_{\lambda}$ (cf \cite{M}). 

On consid\`ere la fonction g\'en\'eratrice

$$\psi(x,y)=\sum_n\frac{1}{n}\sum_{\lambda,\mu\vdash n}
p_{\lambda}(x)p_{\mu}(y)c^n_{\lambda\mu}$$

D'apr\`es (\ref{01}) elle est donn\'ee par

$$\psi(x,y)=\sum_n\sum_{r=0}^{n-1}(-1)^rr!(n-1-r)!\sum_{\lambda,\mu\vdash n}
\frac{p_{\lambda}(x)p_{\mu}(y)}{
z_{\lambda}z_{\mu}}
\chi_{\lambda}^{1^r(n-r)}\chi_{\mu}^{1^r(n-r)}.$$
D'apr\`es \cite{M} I. (7.6) et I.3 example 9, on a 
$$\psi(x,y)=\sum_n\sum_{r=0}^{n-1}(-1)^rr!(n-1-r)!
s_{(n-r-1|r)}(x)s_{(n-r-1|r)}(y)$$
Les $s_{\lambda}$ sont les fonctions de Schur, et $(a|b)=(a+1,1^b)$ suivant la
notation de Frobenius.
D'apr\`es \cite{M} I.3. example 14, on a 
$$\prod_i\frac{1+vx_i}{1-ux_i}=1+(u+v)\sum_{a,b\geq 0}s_{(a,b)}u^av^b$$
donc, en utilisant 
$$\frac{1}{\pi}\int_{\mathbb C}u^k\bar u^le^{-|u|^2}du=\delta_{kl}k!$$ 
on obtient 
\begin{eqnarray*}
\psi(x,y)&=&\frac{1}{\pi^2}\int_{\mathbb C}\int_{\mathbb C}
\left(\frac{\prod_i\frac{1-vx_i}{1-ux_i}-1}{u-v}\right)
\left(\frac{\prod_i\frac{1+\bar vy_i}{1-\bar uy_i}-1}{ \bar u-\bar
v}
\right)e^{-|u|^2-|v|^2}dudv\\
&=&\frac{1}{\pi^2}\int_{\mathbb C}\int_{\mathbb C}
\left(\frac{\exp(\sum_{r}\frac{u^r-v^r}{r}p_r(x))-1}{ u-v}\right)\times\\
&&\qquad
\left(\frac{\exp(\sum_{r}\frac{\bar u^r-(-\bar v)^r}{r}p_r(y))-1}{\bar u+\bar
v}\right)e^{-|u|^2-|v|^2}dudv
\end{eqnarray*}

(l'int\'egrale ne converge pas, mais le d\'eveloppement en s\'erie des
$p_{\lambda}$ converge terme \`a terme).
Faisons le changement de variables $u=a+b, v=a-b$, on trouve
\begin{eqnarray*}
\psi(x,y)&=&\frac{2}{\pi^2}\int_{\mathbb C}\int_{\mathbb C}
\left(\frac{\exp(\sum_{r}\frac{(a+b)^r-(a-b)^r}{r}p_r(x))-1}{ 2b}\right)\times
\\ &&\qquad
\left(\frac{\exp(\sum_{r}\frac{(\bar a+\bar b)^r-(\bar a-\bar b)^r}{
r}p_r(y))-1}{2\bar a}\right)e^{-2|a|^2-2|b|^2}dadb
\end{eqnarray*}
Or le polyn\^ome $Q(a,b)=(a+b)^r-(a-b)^r$ 
a tous ses coefficients positifs, par
cons\'equent quand on d\'eveloppe cette expression en termes des
$p_{\lambda}(x)p_{\mu}(y)$ on trouve des coefficients positifs.
Plus pr\'ecis\'ement, si on note $$R_{\lambda}(a,b)=\frac{1}{b}
\prod_i\frac{
Q^{\alpha_i}_{i}(a,b)}{i^{\alpha_i}\alpha_i!}=\frac{1}{z_{\lambda}b}
\prod_iQ_{\lambda_i}(a,b)$$
pour $\lambda=1^{\alpha_1}2^{\alpha_2}\ldots$, qui est un polyn\^ome homog\`ene 
de degr\'e
$n-1$,  alors on a 
$$c^n_{\lambda,\mu}=\frac{n2^{-n-1}}{\pi^2}
\int_{\mathbb C}\int_{\mathbb C}
R_{\lambda}(a,b)R_{\mu}(\bar b,\bar
a)e^{-|a|^2-|b|^2}dadb$$
ou encore, en appelant $r_{\lambda}^{kl}$ le coefficient de $a^kb^l$ dans
$R_{\lambda}$, (qui est positif) on obtient l'expression
$$c^n_{\lambda,\mu}
=n2^{-n-1}\sum_{k,l}r_{\lambda}^{kl}r_{\mu}^{lk}k!l!$$
 qui est \'equivalente \`a la formule de Goupil et Schaeffer.
\bibliographystyle{amsalpha}

\end{document}